\documentclass[12pt]{amsart}
\usepackage{amssymb,latexsym, amsmath, amscd, array, graphicx}
\usepackage{tikz}
\usetikzlibrary{matrix,arrows}
\usetikzlibrary{positioning}

\swapnumbers \numberwithin{equation}{section}

\theoremstyle{plain}

\newtheorem{thm}{Theorem}[section]
\newtheorem{theorem}[thm]{Theorem}

\newtheorem{lemma}[thm]{Lemma}

\newtheorem{prop}[thm]{Proposition}

\theoremstyle{definition}
\newtheorem{defin}[thm]{Definition}

\DeclareMathOperator{\cat}{{\mbox{\rm cat$_{\rm LS}$}}}

\DeclareMathOperator{\swgt}{\mbox{\rm swgt}}

\DeclareMathOperator{\cuplength}{{\rm cup-length}}

\def\id{\protect\operatorname{id}}

\def\cat{\protect\operatorname{cat}}

\def\Z{{\mathbb Z}}

\def\N{{\mathbb N}}

\def\1{\hbox{\rm\rlap {1}\hskip.03in{\rom I}}}
\def\Bbbone{{\rm1\mathchoice{\kern-0.25em}{\kern-0.25em}
{\kern-0.2em}{\kern-0.2em}I}}

\long\def\forget#1\forgotten{} %

\newcommand\ver[1]{\marginpar{\tiny Changed in Ver \VER}}

\date{\today}
\begin{document}

\title[The LS-category of metric spaces]{The Lusternik-Schnirelmann category of metric spaces}

\author[Tulsi Srinivasan]{Tulsi Srinivasan} %

\address{Tulsi Srinivasan, Department of Mathematics, University
of Florida, 358 Little Hall, Gainesville, FL 32611-8105, USA}
\email{tsrinivasan@ufl.edu}


\begin{abstract}
We extend the theory of the Lusternik-Schnirelmann category to general metric spaces 
by means of covers by arbitrary subsets. We also generalize the definition of the 
strict category weight. We show that if the Bockstein homomorphism on a metric space 
$X$ is non-zero, then $\cat X \geq 2$, and use this to compute the category of Pontryagin 
surfaces.  Additionally, we prove that a Polish space with LS-category $n$ can be presented 
as the inverse limit of ANR spaces of category at most $n$.
\end{abstract}

\maketitle

\maketitle

\section{Introduction}
The {\em Lusternik-Schnirelmann category} (LS-category) of a space $X$ is the smallest integer $n$ such that $X$ admits a cover $U_0,\dots, U_n$ by open sets that are contractible in $X$. It is known that if the sets $U_i$ in the definition are taken to be closed, then the category remains the same for ANR spaces.
For general spaces, the definition of the LS-category by means of closed or open sets is not satisfactory. For instance, the category of the Menger curve is infinite if we consider covers of open or closed subsets. In~\cite{Sr} we considered arbitrary subsets of $X$ in the definition of the LS-category, showed that the new definition coincides with the classical definition on compact metric ANRs, and applied it to fractal Peano continua like Menger spaces and Pontryagin surfaces. In this paper we discuss the LS-category of general metric spaces. In particular, we show that the definition with general subsets agrees with the classical one on all metric ANR spaces. Thus, we may use the notation $\cat X$ for the LS-category defined by means of arbitrary subsets. We note here that our general sets can be assumed to be $G_{\delta}$ on all complete metric spaces.

We then apply the general subsets definition to the category of maps. We recall that the \textit{Lusternik-Schnirelmann category} (LS-category), $\cat f$, of a map $f: X \rightarrow Y$ is the smallest integer $n$ such that $X = \bigcup_{i=0}^n U_i$, where the $U_i$ are open in $X$ and each restriction $f|_{U_i}$ is nullhomotopic. Then the LS-category of a space $X$ can be defined to be $\cat \id_X$. Also, if $A \subset X$, then the relative category of $A$ is $\cat_X A = \cat \iota_A$, where $\iota: A \hookrightarrow X$ is the inclusion map. We generalize this definition by replacing open sets by arbitrary ones (see Definition ~\ref{maps}), and study its properties. We use this to extend the theory of category weight to general metric spaces. We also prove an analogue to a generalization of Freudenthal's theorem in dimension theory, namely that a Polish metric space with LS-category $n$ can be presented as the inverse limit of ANRs of category at most $n$.

Additionally, we apply the theory to Bockstein elements to give a lower bound for the LS-category of Pontryagin surfaces $\Pi_p$, which leads to the equality $\cat \Pi_p=2$ for all $p$. Note that the case $p=2$ was dealt with in~\cite{Sr}, by using cup-length as a lower bound.

\section{Preliminaries}

Recall that given a weakly hereditary topological class of spaces $\mathcal C$, a space $X \in  \mathcal C$ is called an \textit{absolute neighborhood retract} (ANR) \textit{for $\mathcal C$} if every homeomorphic image of $X$ as a closed subset of any $Z \in \mathcal C$ is a neighborhood retract of $Z$. If $X$ is an ANR for the class of metrizable spaces, then we will call $X$ a metric/metrizable ANR.

The following result is known.

\begin{prop}\label{nerve}
Let $X$ be a metric space, $A$ a subset of $X$, and
$\mathcal V' = \{V'_i\}_{i\in \N}$ a cover of $A$ by sets open in $A$ with diam$(V'_i) < \delta_i$. Then
$\mathcal V'$ can be extended to a cover $\mathcal V = \{V_i\}_{i\in N}$ of
$A$ by sets open in $X$ with the same nerve and such that, for all $i\in \N$, $V_i\cap A=V_i'$ and diam$(V_i) < 2\delta_i$.
\begin{proof} 
Set $$V_i = \bigcup_{a \in V'_i}B(a, \frac{1}{2}\min \{\delta_i,d(a, A - V'_i)\}).$$ 
\end{proof}

\end{prop}

\begin{theorem}[Eilenberg-Kuratowski-Wojdyslawski Theorem]\label{KW} ~\cite{B,Hu} Every bounded metric space can be embedded as a closed subset of a convex subset of a Banach space.
\end{theorem}

The following is a version of a lemma that appears in ~\cite{Wa}.

\begin{thm}[Modified Walsh Lemma] \label{Walsh}
Let $X$ be a metric space, $A$ a subset of $X$, $K$ a metrizable ANR,
and $f: A \rightarrow K$ a map. Then, for any $\epsilon: K \rightarrow (0,\infty)$ (not necessarily continuous), there
is an open set $U \supset A$ and a map $g: U \rightarrow K$ such
that:
\begin{enumerate}
    \item For every $u \in U$, there exists $a_u \in A$ such that $d(g(u), f(a_u)) \leq \epsilon(f(a_u))$
    \item $g|_A$ is homotopic to $f$.
\end{enumerate}

\begin{proof} By Theorem ~\ref{KW} (replacing the metric by an equivalent bounded metric if necessary), we can embed $K$ as the closed subset of a convex subset $T$ of a Banach space $Y$. Since $K$ is an ANR,
there is an open neighborhood $O$ of $K$ in $T$, and a retraction
$r:O \rightarrow K$. For $y \in K$, we use the notation $B_T(y,\eta)$ to denote the convex open subset $B(y,\eta)\cap T$ of $T$.

For every $y \in K$,
pick $\delta_y > 0$ such that:
\begin{enumerate}
    \item [(i)]$B_T(y,2\delta_y) \subset O$,
    \item [(ii)] For all $y_1,y_2 \in B_T(y,2\delta_y)$, we have $d(r(y_1),r(y_2)) < \epsilon(y)$ .
\end{enumerate}

For each $a \in A$, pick a neighborhood $U_a$ of $a$ that is open in $A$ so that
$f(U_a) \subset B_T(f(a),\delta_{f(a)})$. Let $\mathcal{V'} =
\{V'_i\}$ be a locally finite refinement of the collection of $U_a$, and
$\mathcal{V} = \{V_i\}$ the cover obtained by applying Proposition ~\ref{nerve} to $\mathcal V'$. Let $U = \bigcup_iV_i$. Note here that we can take the $U_a$ to have as small a diameter as necessary, so the diameter 
of the $V_i$ can be made as small as required. This will be made use of in Lemma ~\ref{extn}.

For each $i$, fix $a_i \in V'_i$ such that $V'_i \subset U_{a_i}$. Let $\{f_{i}: V_i \in \mathcal{V}\}$ be a partition 
of unity dominated by $\mathcal{V}$. Define $h: U \rightarrow Y$ by $h(u) = \sum_{i}f_i(u)f(a_i)$.

Choose any $u \in U$. Then $u$ lies in precisely $k$ elements of $\mathcal V$,
say in $V_{i_1},...,V_{i_k}$. Assume that $\delta_{f(a_{i_j})} =
\max\{\delta_{f(a_{i_1})},...,\delta_{f(a_{i_k})}\}$. Then Proposition ~\ref{nerve} ensures that
$f(a_{i_l}) \in B_T(f(a_{i_j}),2\delta_{f(a_{i_j})})$ for all $l$, so $h(u)$
lies in this convex set as well. This means that $h$ is a map from $U$ to $O$.

Define $g: U \rightarrow K$ by $g = r \circ h$. If $u \in U$, then
we have seen that $h(u)$ lies in $B_T(f(a_u),2\delta_{f(a_u)})$ for
some $a_u \in A$, so $d(g(u),f(a_u)) < \epsilon(f(a_u))$. 

For every $a \in A$, $h(a)$ and $f(a)$ lie in some convex set, so
$h|_A$ is homotopic to $f$ in $O$ via the straight line homotopy.
The composition of this homotopy with $r$ is then a homotopy in $K$
between $g|_A$ and $f$.

\end{proof}
\end{thm}

If $Y$ is a metric ANR with a given open cover $\mathcal U = \{U_\lambda: \lambda \in \Lambda\}$, then there exists a refinement $\mathcal V$ of $\mathcal U$ such that given any $X$ and any $f,g: X \rightarrow Y$ such that $f$ and $g$ are $\mathcal V$-close, the maps $f$ and $g$ are $\mathcal U$-homotopic ~\cite{Hu}. We require a far less general version of this fact:

\vspace{3mm}

\begin{thm} ~\cite{Hu}
\label{epsilon} Let $K$ be a metric ANR. Then there exists a (not necessarily continuous) map $\epsilon: K \rightarrow (0,\infty)$ such that the covering $\mathcal W = \{B(y, \epsilon(y)): y \in K\}$ has the following property: for any metric space $X$, any $\mathcal W$-close maps $f,g: X \rightarrow K$ are homotopic to each other.
\end{thm}

Let $K$ be a metric ANR, and $PK$ the path space $\{\phi: [0,1] \rightarrow K|\phi(1) = y_{0}\}$ for some $y_{0} \in K$, equipped with the sup norm. For any $\epsilon: K \rightarrow(0,\infty)$, define $\bar \epsilon: PK \rightarrow (0,\infty)$ by $\bar\epsilon(\phi) = \epsilon(\phi(0))$. A metrizable space is an ANR for metrizable spaces iff it is an absolute neighborhood extensor (ANE) for metrizable spaces ~\cite{Hu}, from which it follows quite easily that if $K$ is a metric ANR, then so is $PK$. 

\begin{lemma}\label{extn}
Let $X$ be a metric space, $A \subset X$, $K$ a metric ANR, and $f: X \rightarrow K$ a map such that the restriction of $f$ to $A$ is nullhomotopic. Then there exists $U \supset A$ open
in $X$ such that the restriction of $f$ to $U$ is also nullhomotopic.

\begin{proof} Let $D$ be the metric on $K$ and $d$ the metric on $X$. Let $PK$ be the path space
$\{\phi: [0,1] \rightarrow K| \phi(1) = y_{0}\}$, under the sup metric $D'$. Let $\epsilon: K \rightarrow (0,\infty)$ be a map as in Theorem ~\ref{epsilon} and $\bar \epsilon: PK \rightarrow (0,\infty)$ as defined above.

Since $f|_{A}$ is nullhomotopic, there is a map $F: A \rightarrow PK$ satisfying $F(a)(0) = f(a)$ for all 
$a \in A$. For each $a \in A$, pick $\delta_a >0$ so that for all $x \in X$ with $d(x,a)< \delta_a$, we have  $D(f(x),f(a)) < \epsilon(f(a))$. We apply Theorem ~\ref{Walsh} to $F$, and construct an open neighborhood $U$ of $A$ and a map $G: U \rightarrow PK$ such that for every $u \in U$, there exists $a_u \in A$ such that
$D'(G(u),F(a_u)) < \bar \epsilon (F(a_u))$. As noted in the proof of Theorem ~\ref{Walsh}, we can construct $U$ so that diam$(V_i) < \delta_{a_i}$ for all $i$, the $a_i$ and $V_i$ being as in the proof of Theorem ~\ref{Walsh}, and $a_u$ being one of the $a_i$.

Let $g: U \rightarrow K$ be given by $g(u) = G(u)(0)$. Then $D(g(u),f(a_u)) \leq \bar \epsilon (F(a_u)) = \epsilon(f(a_u))$ and $D(f(u), f(a_u)) <  \epsilon(f(a_u))$, and so both $g(u)$ and $f(u)$ lie in $B(f(a_u),\epsilon(f(a_u)))$. It follows that $f|_U$ is homotopic to $g$, which, by definition, is nullhomotopic.
\end{proof}
\end{lemma}

\section{Generalized category of maps}

\begin{defin}\label{maps}
  For a map $f: X \rightarrow Y$, define the
{\em generalized LS-category} of $f$, $\cat_g f$, to
be the smallest integer $k$ such that $X = \bigcup_{i=0}^{k}A_{i}$,
where each restriction $f|_{A_{i}}$ is nullhomotopic. If no such $k$ exists, we 
say that the category is infinite.
\end{defin}

By definition, we have $\cat_g X=\cat_g \id_X$.

\begin{thm} If $f: X \rightarrow Y$ be a map between a metric space $X$ and a metric ANR $Y$, then 
$\cat_gf = \cat f$. In particular, if $X$ is a metric ANR, then $\cat_g X = \cat X$.
\begin{proof} Clearly $\cat_gf \leq \cat f$, and equality holds if $\cat_gf$ is infinite. Suppose 
$\cat_gf = n$. Then we can write $X = \bigcup_{i = 0}^{n}A_i$, where each $f|_{A_i}$ is 
nullhomotopic. By Lemma~\ref{extn}, there exist open sets $U_i$ containing $A_i$ such that 
each $f|_{U_i}$ is nullhomotopic. It follows that $\cat f \leq n.$
\end{proof}
\end{thm}

In the light of this, we will use the notation $\cat X$ for $\cat_gX$ for the remainder of this paper. However, 
we continue to use the subscript when referring to the LS-category of maps. We now list some basic properties of $\cat_g$. The proofs are omitted, as they are nearly identical to those for $\cat$ (see, for instance, ~\cite{CLOT}).

\begin{prop} [Properties of $\cat_g $] \label{properties}

i) For any map $f: X \rightarrow Y$, we have $\cat_g f \leq \min\{\cat X,\cat Y\}$\\

ii) If $f: X \rightarrow Y$ and $h: Y \rightarrow Z$ are maps, then $$\cat_g (h\circ f) \leq \min\{\cat_g f, \cat_g h\}$$

iii) If $f: X \rightarrow Y$ is a homotopy equivalence, then $\cat_g f = \cat X = \cat Y$\\

iv) If $F \overset{i}{\rightarrow} E \overset{p}{\rightarrow} B$ is a fibration, then

$$\cat(E) \leq (\cat_g(i) + 1)(\cat_g(p) + 1) - 1.$$

\end{prop}

Given a commutative ring $R$, the $R$-cup-length, $\cuplength_R X$, of a space $X$ is the smallest integer $k$ 
such that all cup-products of length $k+1$ vanish in the cohomology ring $H^*(X;R)$. We show
that cup-length is a lower bound for $\cat_g$.

\vspace{3mm}

\begin{theorem}~\cite[Theorem I.4.1 and Theorem 2.2.8, Appendix 1]{MS} Every CW complex has the homotopy type of a metric ANR.\label{CW}
\end{theorem}

\begin{theorem} Let $f: X \rightarrow Y$ be a map between metric spaces. Then for any commutative ring $R$, $\cat_g(f) \geq \cuplength_R(Im(f^*))$. In particular, $\cat X \geq \cuplength_R X$.

\begin{proof} We consider $X$ and $Y$ with the Alexander-Spanier cohomology or \v Cech cohomology (which coincide on paracompact spaces ~\cite{HDD}). Let $\cat_g(f) = n$, and write $X = \cup_{i=0}^n A_i$, where each restriction $f|_{A_i}$ is nullhomotopic. Suppose $\cuplength_R(Im(f^*)) > n$. Then there exist $\alpha_i \in H^{k_i}(Y;R)$ such that $f^*(\alpha_0) \smile ... \smile f^*(\alpha_n) \neq 0$. By the representability of the \v Cech cohomology,  there exist maps $g_i: Y \rightarrow K(R,k_i)$ such that $g^*(u) = \alpha_i$, for a certain $u \in H^{k_i}(K(R,k_i);R)$. By Theorem ~\ref{CW} we may assume that the $K(R;k_i)$ are metrizable. Then each map $g_i \circ f: X \rightarrow K_i$ is a map between a metric space and a metric ANR with $g_i \circ f|_{A_i}$ nullhomotopic, so Lemma ~\ref{extn} implies that there exist open sets $U_i \supset A_i$ such that $g_i\circ f|_{U_i}$ is nullhomotopic. Let $j_i: U_i \rightarrow X$ be the inclusion, and consider the following exact sequence:

\begin{tikzpicture}[align=center, node distance=3cm, auto]
  \node (A) {$...$};
  \node (B) [right of=A, xshift=-1cm]{$H^{k_i}(X,U_i;R)$};
  \node (C) [right of=B] {$H^{k_i}(X;R)$};
  \node (D) [right of=C] {$H^{k_i}(U_i;R)$};
  \node (E) [right of=D, xshift=-1cm] {$...$};
  \node (F) [below of=C, yshift=1.5cm] {$H^{k_i}(Y;R)$};
  \draw[->] (A) to node {$ $} (B);
  \draw[->] (B) to node {$q_i^*$} (C);
  \draw[->] (C) to node{$j_i^*$} (D);
  \draw[->] (D) to node {$ $} (E);
  \draw[->] (F) to node {$f^*$} (C);
\end{tikzpicture}

We have $j_i^*(f^*(\alpha_i)) = (g_i \circ f \circ j)^*(u) = 0$, so exactness implies that there exist $\bar \alpha_i \in H^{k_i}(X,U_i;R)$ such that $q_i^*(\bar \alpha_i) = f^*(\alpha_i)$. The rest of the proof follows as in the proof of the analogous statement for $\cat$ (see ~\cite{CLOT}). Namely, in view of the relative 
cup-product formula for the Alexander-Spanier cohomology (see ~\cite{Sp}, p. 315), 
$$H^k(X,U;R)\times H^l(X,V;R)\stackrel{\cup}{\longrightarrow} H^{k+l}(X,U\cup V;R),$$
and the fact that $H^n(X,X;R)=0$, we obtain a contradiction.

\end{proof}
\end{theorem}

We conclude this section with a result analogous to the following generalization of Freudenthal's theorem in dimension theory: an $n$-dimensional completely metrizable space can be expressed as the inverse limit of polyhedra of dimension at most $n$ ~\cite{Na}.

\begin{theorem} Let $X$ be a Polish space with $\cat X \leq n$. Then $X$ can be presented as the inverse limit of a sequence of ANR spaces, each of which has category $\leq n$. 

\begin{proof} Suppose $X = \cup_{i=0}^n A_i$, where each $A_i$ is contractible in $X$. Embed $X$ as a closed subset of $\mathbb R^\omega$ ~\cite{Ch}, and take an open neighborhood $U$ of $X$. By Lemma ~\ref{extn}, we can enlarge each $A_i$ to $U_{1,i}$ that is open and contractible in $U$, with $U_{1,i} \subset N_1(A_i)$.  Set  $U_1 = \bigcup_i U_{1,i}$  and  $L_1 = \bigcup_i CU_{1,i}$, where $CY =  (Y\times [0,1])/(Y \times \{1\})$ denotes the cone on $Y$. Then $U_1$ and $L_1$ are ANR spaces with $\cat L_1 \leq n$.

Next, we expand each $A_i$ to $U_{2,i}$ that is open and contractible in $U_1$, with $U_{2,i} \subset N_{1/2}(A_i)$.  Set $U_2 = \bigcup_i U_{2,i}$ and $L_2 = \bigcup_i CU_{2,i}$. Then $U_2$ and $L_2$ are ANRs with $\cat L_2 \leq n$. Since each $U_{2,i}$ is contractible in $U_1$, there is a map $H_{2,i}: CU_{2,i} \rightarrow U_1$ with $H_{2,i}|_{U_{2,i} \times \{0\}} = \id_{U_{2,i}}$. We paste these together to obtain a map $q_2: L_2 \rightarrow L_1$ with $q_2|_{U_2} = \id_{U_2}$. Continuing this way, we obtain the inverse system $\{L_j,q_j\}$ of ANR spaces of dimension $\leq n$, which we claim has inverse limit $L$ identical to $X$. 

Since $X$ is complete and each $U_{n,i} \subset N_{1/n}(A_i)$, the inverse limit of the inverse sequence $\{U_j,q_j|_{U_j}\}$ is $X$. Thus the maps $q_n: L_n \rightarrow U_{n-1} \subset L_{n-1}$ induce a map $q: L \rightarrow X$. Also, the inclusions $X \hookrightarrow L_n$ induce $p: X \rightarrow L$ with $q \circ  p  = \id_X$, which implies that $q$ is surjective. Finally, suppose two elements $(y_j)$ and $(z_j)$ in $L$ are both mapped by $q$ to $(x_j) \in X$. Then for each $n$ we have $y_n = q_{n+1}(y_{n+1}) = x_n = q_{n+1}(z_{n+1}) = z_n$, so $q$ is injective as well.

\end{proof}
\end{theorem}

\section{Category and the Bockstein homomorphism}

The Bockstein homomorphism on a space $X$ is the map $$\beta_X: H^1(X; \Z_p) \rightarrow H^2(X;\Z_p)$$ obtained from the exact sequence $0 \rightarrow \Z_p \overset{p}{\rightarrow} \Z_{p^2} \rightarrow \Z_p \rightarrow 0$. It is known that for a metric ANR $X$, if $\beta_X \neq 0$, then $\cat X \geq 2$. ~\cite{FH, Ru}

For metric ANR spaces, the usual definition of category coincides with Ganea's definition, i.e., given $f: X \rightarrow Y$, $\cat f \leq n$ iff there exists a map $\tilde f: X \rightarrow G_n(Y)$ such that $p_n \circ \tilde f \simeq f$, where $G_n(Y)$ is the $n$th Ganea space of $Y$ and $p_n: G_n(Y) \rightarrow Y$ the $n$th Ganea map ~\cite{CLOT}. The Ganea definition of category, together with the fact that $G_n(K(\Pi,1))$ has the homotopy type of a CW-complex of dimension $\leq n$ for any group $\Pi$, yields a proof that $\beta_X \neq 0 \Rightarrow \cat X \geq 2$. It is not clear whether or not the Ganea definition can be applied to general spaces, but we can modify the technique above to obtain the same lower bound on the LS-category of a general metric space.

\begin{theorem} Let $X$ be a metrizable space on which the Bockstein homomorphism $\beta_X$ is non-zero. Then $\cat X \geq 2$.\label{Bockstein}

\begin{proof} Since $\beta_X: H^1(X;\Z_p) \rightarrow H^2(X;\Z_p) \neq 0$, there exists $\alpha \in  H^1(X;\Z_p)$ such that $\beta_X(\alpha) \neq 0$. By the representability of the \v Cech cohomology, there is a map $f: X \rightarrow K(\Z_p,1)$ such that $f^*(i) = \alpha$, for a certain $i \in H^1(K(\Z_p,1);\Z_p)$.

By Theorem ~\ref{KW}, we can embed $X$ as a closed subset of a convex subset $T$ of a Banach space. Since $K(\Z_p,1)$ is an ANE, $f$ extends to a map $\tilde f: L \rightarrow K(\Z_p,1)$ over some open neighborhood $L$ of $X$ in $T$. Since $T$ is a convex subset of a Banach space, the Dugundji extension theorem ~\cite{Du} implies that it is an ANR (in fact, an AR), and therefore $L$ is also an ANR ~\cite{Hu}. 

Suppose, for a contradiction, that $\cat X < 2$. Then $X = A_1 \cup A_2$, where the $A_i$ are contractible in $X$. By Lemma ~\ref{extn}, the $A_i$ can be expanded to open sets $U_i$ that are contractible in $L$. Set $K = U_1 \cup U_2$. Let $\iota_K: K \rightarrow L$ and $\iota_X: X \rightarrow K$ be the inclusion maps, and 
let $\bar f: K \rightarrow K(\Z_p,1)$ be the restriction of $\tilde f$ to $K$. We claim that the Bockstein homomorphism $\beta_K: H^1(K;\Z_p) \rightarrow H^2(K;\Z_p)$ is non-zero. 

Let $\bar \alpha = \bar f^*(i)$. From the diagram below, we have $0 \neq \beta_X(\alpha) = \beta_X \circ f^*(i) = \beta_X \circ \iota_X^* \circ \bar f^*(i) = \iota_X^*\circ \beta_K \circ \bar f^*(i) = \iota_X^* (\beta_K(\bar \alpha))$.

\begin{equation*}
\begin{CD}
H^1(K;\Z_p) @> \iota_X^* >> H^1(X;\Z_p)\\
@VV \beta_K V @VV \beta_X V \\
H^2(K;\Z_p) @> \iota_X^* >> H^2(X;\Z_p)\\
\end{CD}
\end{equation*}

This proves the claim.

However, since $\cat \iota_K < 2$ and $K,L$ are metric ANRs, there exists a section $s$ from $K$ to the Ganea space  $G_1(L)$, which induces $s^*:H^2(G_1(L);\Z_p) \rightarrow  H^2(K;\Z_p)$. Since $G_1(K(\Z_p,1))$ has the homotopy type of a cell-complex of dimension $\leq 1$, $H^2(G_1(K(\Z_p,1));\Z_p) = 0$, so $i$ is mapped to $0$ in $H^2(K;\Z_p)$ under $s^* \circ G_1(\tilde f^*) \circ p_1^*\circ \beta$ (see the diagram below). This implies that $\beta_K(\bar \alpha) = 0$, which is a contradiction.

\begin{tikzpicture}[align=center, node distance=4.5cm, auto]
  \node (A) {$H^2(G_1(K(\Z_p,1));\Z_p)$};
  \node (B) [right of=A] {$H^2(G_1(L);\Z_p)$};
  \node (C) [right of=B] {$H^2(G_1(K);\Z_p)$};
  \node (D) [below of=A, yshift=3cm] {$H^2(K(\Z_p,1);\Z_p)$};
  \node (E) [below of=B, yshift=3cm] {$H^2(L;\Z_p)$};
  \node (F) [below of=C, yshift=3cm] {$ H^2(K;\Z_p)$};
  \node (G) [below of=D, yshift=3cm] {$H^1(K(\Z_p,1);\Z_p)$};
  \node (H) [below of=E, yshift=3cm] {$H^1(L;\Z_p)$};
  \node (I) [below of=F, yshift=3cm] {$H^1(K;\Z_p)$};
  \draw[->] (A)to node {$G_1(\tilde f^*)$} (B);
  \draw[->] (B)to node{$G_1(\iota_K^*)$} (C);
  \draw[->] (D)to node {$\tilde f^*$} (E);
  \draw[->] (E) to node {$\iota_K^*$} (F);
  \draw[->] (G)to node {$\tilde f^*$} (H);
  \draw[->] (H) to node {$\iota_K^*$} (I);
  \draw[<-] (A) to node {$p_1^*$} (D);
  \draw[<-] (B) to node {$p_{1,L}^*$} (E);
  \draw[<-] (C) to node {$p_{1,K}^*$} (F);
  \draw[<-] (D) to node {$\beta$} (G);
  \draw[<-] (E) to node {$\beta_L$} (H);
  \draw[<-] (F) to node {$\beta_K$} (I);
  \draw[->, dashed] (B) to node {$s^*$} (F);

\end{tikzpicture}

\end{proof}
\end{theorem}

We use this result to compute the category of the Pontryagin surface $\Pi_p$, which is constructed as follows: let $M_p$ be the mapping cylinder of the degree $p$ map $f_p: S^1 \rightarrow S^1$, and let $\partial M_p \subset M_p$ denote the domain of $f_p$.  Let $L_1 = S^2$, and suppose $L_{k-1}$ has been constructed for $k \geq 2$. To construct $L_k$, take a sufficiently small triangulation of $L_{k-1}$, remove the interior of each of the triangles $\Delta$ and identify each $\partial \Delta$ with a copy of $\partial M_p$. Define $q_{k+1}: L_{k+1} \rightarrow L_k$ to be the map that collapses each $image(f_p) \subset M_p$ to a single point. $\Pi_p$ is defined to be the inverse limit of $\{L_k, q_{k+1}\}$.

\begin{theorem} The category of the Pontryagin surface $\Pi_p = 2$.

\begin{proof} We refine the construction above by attaching copies of $M_p$ one at a time as follows: Let $M_p$ be as above, and let $X_1 = L_1 = S_2$. To construct $L_2$, we triangulate $S_2$ with, say, $n_1$ triangles. Let $X_2$ be the space obtained when one of the triangles $\Delta$ has its interior removed and its boundary identified with $\partial M_p$. Define $\lambda_2: X_2 \rightarrow X_1$ to be the map that collapses $image(f_p) \subset M_p$ to a point. Let $X_3$ be obtained from $X_2$ by removing the interior of a second triangle and attaching a second copy of $M_p$, and so on. We get $X_{n_1} = L_2$, and then use the same process to construct $X_{n_1 + 1}$. Let $\lambda_{k+1}: X_{k+1} \rightarrow X_k$ be the map that collapses $image(f_p)$ in the copy of $M_p$ added to $X_k$ to form $X_{k+1}$ to a single point. Then $\Pi_p$ is the inverse limit of $\{X_k,\lambda_{k+1}\}$.

We claim that the Bockstein homomorphism on $\beta$ on $\Pi_p$ is non-zero. This, along with the fact that $\cat X \leq \dim X$ ~\cite{Sr}, will imply that $\cat \Pi_p = 2$. 

The Mayer-Vietoris sequence can be used to show that the maps $\lambda_{k+1}^*: H^2(X_k;\Z_p) \rightarrow H^2(X_{k+1};\Z_p)$ are isomorphisms on $\Z_p$ for all $k$ ~\cite{Dr}, and the long exact sequence for cohomology implies that each $\lambda_{k+1}^*: H^1(X_k;\Z_p) \rightarrow H^1(X_{k+1};\Z_p)$ is injective. Also, since $H^2(X_2)$ has one $\Z_p$ component, the Bockstein map $\beta_2: H^1(X_2;\Z_p) \rightarrow H^2(X_2;\Z_p)$ is an isomorphism ~\cite{H}, so the commutative diagram below can be used to show, inductively, that each $\beta_k: H^1(X_k;\Z_p) \rightarrow H^2(X_k;\Z_p)$ is non-zero:

\begin{tikzpicture}[align=center, node distance=4cm, auto]
  \node (A) {$H^1(X_n;\Z_p)$};
  \node (B) [right of=A]{$H^2(X_n;\Z_p)$};
  \node (C) [below of =A, yshift=2cm] {$H^1(X_{n+1};\Z_p)$};
  \node (D) [below of =B, yshift=2cm]{$H^2(X_{n+1};\Z_p)$};
  \draw[->] (A) to node {$\beta_n$} (B);
  \draw[->] (C) to node {$\beta_{n+1}$} (D);
  \draw[->] (A) to node {$ \lambda_{n+1}^*$} (C);
  \draw[->] (B) to node {$ \cong $} (D);
 
\end{tikzpicture}

\end{proof}
\end{theorem}

\section{Category Weight} For a metric space $X$ the \emph{strict category weight}, $\swgt(u)$, of $0 \neq u \in H^*(X;R)$, is the largest integer $k$ such that for every map $\phi: A \rightarrow X$, where $A$ is a metric space and $\cat \phi < k$, we have $\phi^*(u) = 0$. The strict category weight was first introduced in ~\cite{Ru}, in which $X$ and $A$ are assumed to be CW complexes, and in which it is shown that $\cat X \geq \swgt(u)$ for every $0 \neq u \in H^*(X;R)$. Strict category weight is a homotopy invariant version of Fadell and Husseini's category weight ~\cite{FH}, which is also a lower bound for the LS-category of a space. The definition can also be given with $X$ and $A$  paracompact Hausdorff spaces, as in ~\cite{Ru1}. We generalize this as follows:

\begin{defin} \label{wgt} Let $X$ be a metric space and let $0 \neq u \in H^*(X;R)$. Define the \em{generalized strict category weight} of $u$, $\swgt_g(u)$, by 
$$\swgt_g(u) \geq k \Leftrightarrow \phi^*(u) = 0 \mbox{ for all } \phi: A \rightarrow X \mbox{ with }  A \mbox{ metric,  }   \cat_g\phi <k. $$
\end{defin}

Here $R$ is a commutative ring, and $H^*(X;R)$ is the Alexander-Spanier or C\v ech cohomology.

\begin{prop}\label{comparision} Let $X$ be a metric space and $0 \neq u \in H^*(X;R)$. Then 

i) If $f: Y \rightarrow X$ is such that $f^*(u) \neq 0$, then $\swgt_g(u) \leq \cat_g f$. In particular, $\swgt_g(u) \leq \cat X$.\\
ii) If $X$ is an ANR, then $\swgt_g(u) = \swgt(u)$.

\begin{proof}
i) This follows from definition.\\
ii) Clearly, $\swgt_g(u) \leq \swgt(u)$. Suppose $X$ is an ANR and $\swgt(u) = k$. Let $\phi: A \rightarrow X$ be a map from a metric space $A$ with $\cat_g \phi < k$. By Lemma ~\ref{extn}, $\cat \phi < k$, and hence $\phi^*(u) = 0$. 
\end{proof}
\end{prop}

\begin{prop} Let $u,v \in H^*(X;R)$. Then $$\swgt_g(u \smile v) \geq \swgt_g(u) + \swgt_g(v).$$
\begin{proof} Let $\swgt_g(u) = m$ and $\swgt_g(v) = n$, and consider a map $\phi: A \rightarrow X$ with $\cat_g\phi < m + n$. We need to show that $\phi^*(u \smile v) = 0$. 

Let $u = f^*(u')$ and $v = g^*(v')$ for $f: X \rightarrow K(R,k_1)$ and $g: X \rightarrow K(R,k_2)$, and suppose $A = B_1 \cup ...\cup B_n \cup C_1 \cup...\cup C_m$, where the restriction of $\phi$ to each $B_i$ and $C_j$ is nullhomotopic. Let $B = \bigcup_{i=1}^nB_i$ and $C = \bigcup_{j=1}^mC_j$. Since $\cat_g(\phi|_B) < n$, we have $(\phi|_B)^*(u) = 0$, and hence the restriction of $f\circ \phi$ to $B$ is nullhomotopic. Then Lemma ~\ref{extn} gives us an open set $U \supset B$, the restriction of $f\circ\phi$ to which is nullhomotopic. Similarly, we obtain $V \supset C$ such that the restriction of $g\circ\phi$ to $V$ is nullhomotopic.

Consider the exact sequence associated with $(A,U)$:

\begin{tikzpicture}[align=center, node distance=3cm, auto]
  \node (A) {$...$};
  \node (B) [right of=A, xshift=-1cm]{$H^{k_1}(A,U;R)$};
  \node (C) [right of=B] {$H^{k_1}(A;R)$};
  \node (D) [right of=C] {$H^{k_1}(U;R)$};
  \node (E) [right of=D, xshift=-1cm] {$...$};
  \draw[->] (A) to node {$ $} (B);
  \draw[->] (B) to node {$q_1^*$} (C);
  \draw[->] (C) to node{$j^*$} (D);
  \draw[->] (D) to node {$ $} (E);

\end{tikzpicture}

Since $j^*\phi^*(u) = 0$, there exists a $\bar u$ with $q_1^*(\bar u) = \phi^*(u)$. Similarly, there exists a $\bar v \in H^{k_2}(A,V;R)$ with $q_2^*(\bar v) = \phi^*(v)$. Now $\bar u \smile \bar v = 0$, as it lies in $H^*(A,U\cup V;R)$. But $\bar u \smile \bar v \mapsto \phi^*(u) \smile \phi^*(v) = \phi^*(u \smile v)$, proving the claim.

\end{proof}
\end{prop}

\begin{prop} Let $X$ be a metric space. Then\\
i) If $f: Y \rightarrow X$ is such that $f^*(u) \neq 0$ for some $u \in H^*(X;R)$, then $\swgt_g(f^*(u)) \geq \swgt_g(u)$. If $f$ is a homotopy equivalence, then $\swgt_g(f^*(u)) = \swgt_g(u)$.\\
ii) For $0 \neq w \in H^s(K(\Pi,1);R)$, we have $\swgt_g(w) \geq s$.\\
iii) If $f: X \rightarrow K(\Pi,1)$ and $0 \neq f^*(w) = u \in H^s(X;R)$, then  $\swgt_g(u) \geq s$.\\
iv) If $\beta_X(u) \neq 0$ for some $u \in H^1(X;\Z_p)$, then $\swgt_g(\beta_X(u)) \geq 2$.

\begin{proof}
i) Suppose $\swgt_g(u) = n$, and consider $\phi: A \rightarrow Y$ with $\cat_g\phi < n$. By Proposition ~\ref{properties}, $\cat_g (f\circ \phi) < n$, so $0 = (f\circ \phi)^*(u) = \phi^*(f^*(u))$. The second statement follows from the first.

ii) By Theorem ~\ref{CW}, we can pick a metric $K(\Pi,1)$ space. Then $\swgt_g(u) = \swgt(u)$. It is known that if $0 \neq u \in H^s(K(\Pi,1);R)$, then $\swgt(u) \geq s$ ~\cite{CLOT}. 

iii) This follows from i) and ii).

iv) Suppose $u = f^*(i)$, where $f: X \rightarrow K(\Z_p,1)$. Then $\beta_X(u) = f^*(\beta_{K(\Z_p,1)}(i))$, so we can apply iii).
\end{proof}
\end{prop}

Note that iv), together with the inequality $\swgt_g(u) \leq \cat  X$, provides an alternative proof of Theorem ~\ref{Bockstein}.

\section{Acknowledgement}
I am grateful to my advisor, Alexander Dranishnikov, for all his ideas and insights, and also to Yuli Rudyak for suggesting that I study category weight.

\end{document}